\documentclass[12pt]{amsart} 
\usepackage[utf8]{inputenc}
\usepackage{amsmath,amssymb,amsthm,graphicx,multirow}
\usepackage[left=2cm,right=2cm,top=2cm,bottom=2cm]{geometry}

\title{A Reciprocity Formula on Multicurves}

\newtheorem{theorem}{Theorem}[section]
\newtheorem{lemma}[theorem]{Lemma}
\newtheorem{proposition}[theorem]{Proposition}

\theoremstyle{definition}
\newtheorem{definition}[theorem]{Definition}
\newtheorem{example}[theorem]{Example}

\newtheorem*{acknowledgment}{Acknowledgment}

\theoremstyle{remark}
\newtheorem{remark}[theorem]{Remark}

\numberwithin{equation}{section}

\begin{document}
\def\@maketitle{%
  \normalfont\normalsize
  \let\@makefnmark\relax  \let\@thefnmark\relax
\ifx\@empty\@date\else \@footnotetext{\@setdate}\fi
  \ifx\@empty\@subjclass\else \@footnotetext{\@setsubjclass}\fi
  \ifx\@empty\@keywords\else \@footnotetext{\@setkeywords}\fi
  \ifx\@empty\thankses\else \@footnotetext{%
    \def\par{\let\par\@par}\@setthanks}\fi
  \@mkboth{\@nx\shortauthors}{\@nx\shorttitle}%
\global\topskip42\p@ % 5.5 picas to the base of the first title line
  \@settitle
  \ifx\@empty\authors \else \@setauthors \fi
  \ifx\@empty\@dedicatory
  \else
    \baselineskip18\p@
    \vtop{\centering{\footnotesize\itshape\@dedicatory\@@par}%
      \global\dimen@i\prevdepth}\prevdepth\dimen@i
  \fi
  \normalsize
  \dimen@34\p@ \advance\dimen@-\baselineskip
  \vskip\dimen@\relax
} % end \@maketitle

\author{Juhan Kim}
\address{Department of Mathematics, Seoul National University, 1 Gwanak-ro, Gwanak-gu, Seoul 08826}
\email{paytech@snu.ac.kr}

\maketitle

\begin{abstract}
Given a specific collection of curves on an oriented surface with punctures, we associate a power series by counting its intersections with multicurves. This paper presents a reciprocity formula on the power series when multicurves with no component contractible to a puncture are concerned, as a generalization of the reciprocity presented in \cite{MR4157427}.
\end{abstract}

\section{Introduction} 
\begin{figure}[h]
    \centering
    \includegraphics[scale=0.75]{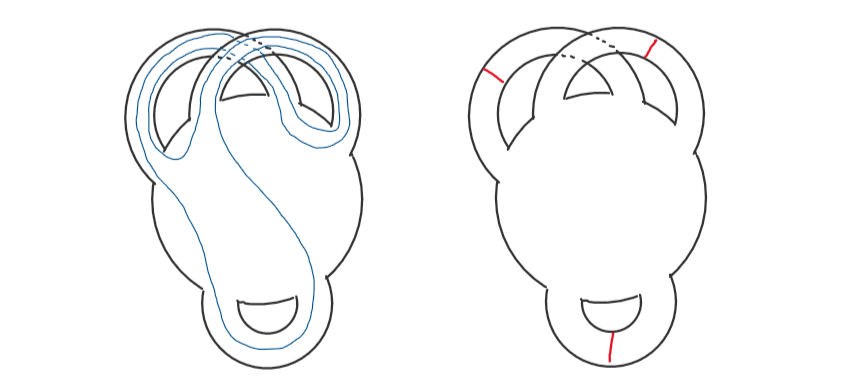}
    \caption{Examples of a multicurve and a counting curve on a ribbon graph}
    \label{fig:ex0}
\end{figure}

The following reciprocity formula Theorem~\ref{thm:0} is originally introduced in \cite{MR4157427}. It is presented as a topological interpretation of log Calabi-Yau property of a relative character variety of a surface. Let there be a closed disk. Attach some ``handles'' to the disk in an orientable way to obtain a ribbon graph $R$, just as in Figure~\ref{fig:ex0}. Define \emph{multicurve} on $R$ to be a finite disjoint union of simple closed non-contractible curves on $R$. We say a multicurve is \emph{essential} when it has no component isotopic to a boundary component of $R$. For $r \in \mathbb{Z}_{\geq 0}$, Define $c(r)$ to be the number of essential multicurves (up to isotopy) that passes through the handles $r$ times in total. More details can be found in \cite{MR4157427}.

\begin{theorem}
\label{thm:0}
The series $Z(t)=\sum\limits_{r=0}^\infty c(r) t^r$is a rational function, and $$Z(t^{-1})=Z\left(t\right).$$
\end{theorem}

The goal of this paper is to provide a new proof of this formula, furthermore generalizing it to a multivariate version on arbitrary surfaces with boundary or punctured surfaces. Let $S_{g,n}$ be an oriented surface of genus $g$ with $n\geq1$ distinct punctures. First we define multicurve in the identical way. 

\begin{definition}[multicurve]
A \emph{multicurve} on $S_{g,n}$ is a finite disjoint union of simple closed non-contractible curves on $S_{g,n}$. A multicurve is \emph{essential} when it has no component contractible to a puncture. We allow a trivial multicurve which is an empty set. Denote by $M$ the set of isotopy classes of multicurves, by $E$ those of essential multicurves. 
\end{definition}

In the case of ribbon graph, we have counted the number of curves passing through handles. Since this amounts to counting intersections of a multicurve and another collection of curves that cut across each handle(Figure~\ref{fig:ex0}), we generalize this to define the following notion of counting curve. 

\begin{definition}[counting curve]
Consider a finite disjoint union of simple curves on  $S_{g,n}$ that connect between the punctures. When such disjoint union meets every nontrivial multicurve, we say that union is a \emph{counting curve}.  
\end{definition}

\begin{example}
In Figure~\ref{fig:ex1}, $a\cup b\cup c$ is a counting curve on $S_{1,2}$. Its components are simple non-contractible curves, and it meets every nontrivial multicurve. $a\cup b\cup d$ is not a counting curve. Even though its component are simple and non-contractible, it does not meet $m$, a nontrivial multicurve. 
\begin{figure}[h]
    \centering
    \includegraphics[scale=0.75]{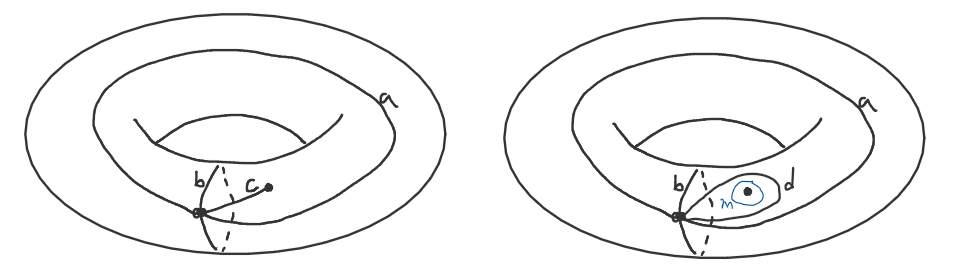}
    \caption{Example and non-example of counting curve on $S_{1,2}$}
    \label{fig:ex1}
\end{figure}
\end{example}

\begin{definition}

Given a counting curve $C$ with $m$ components $C_1,\ldots,C_m$, we define a counting function $p:M\rightarrow \mathbb{Z}^{m}_{\geq 0}$ as $p(x)_i=\min\limits_{c \in x} \# (c\cap C_i) \ (x \in M)$. Let $G_\alpha=\#\left\{x \in M\ |\ p(x)=\alpha \right\}$, $F_\alpha=\#\left\{x \in E\ |\ p(x)=\alpha \right\}$. Now define the following power series. $$g_C(x)=\sum\limits_{\alpha \in \mathbb{Z}^{m}_{\geq0}}G_\alpha x^\alpha=\sum\limits_{\alpha \in \mathbb{Z}^{m}_{\geq0}}G_\alpha x_1^{\alpha_1}x_2^{\alpha_2}\ldots x_m^{\alpha_m}$$ $$f_C(x)=\sum\limits_{\alpha \in \mathbb{Z}^{m}_{\geq0}}F_\alpha x^\alpha=\sum\limits_{\alpha \in \mathbb{Z}^{m}_{\geq0}}F_\alpha x_1^{\alpha_1}x_2^{\alpha_2}\ldots x_m^{\alpha_m}$$

\end{definition}
 Then we are ready to state the main theorem. The fact that $g_C$ and $f_C$ are well-defined (that is, $G_\alpha$ and $F_\alpha$ are always finite) will be evident in the proof of Theorem~\ref{thm:main}.

\begin{theorem}
\label{thm:main}
$g_C$ and $f_C$ are rational functions, and $$f_C(x^{-1})=f_C(x)$$
where $f_C(x^{-1})=f_C(x_1^{-1},x_2^{-1},\ldots,x_m^{-1})$.
\end{theorem}

\begin{example}
A counting curve $a \cup b$ on $S_{1,1}$ is described in Figure~\ref{fig:ex3}. One easily checks that there is a bijection $\phi: \mathbb{Z}^2/\left\{\pm 1\right\} \to E$ such that $p\left((s,t)\right)=\left(\left|s\right|,\left|t\right|\right)$, as suggested in Figure~\ref{fig:ex3}.
\begin{figure}[h]
    \centering
    \includegraphics[scale=0.75]{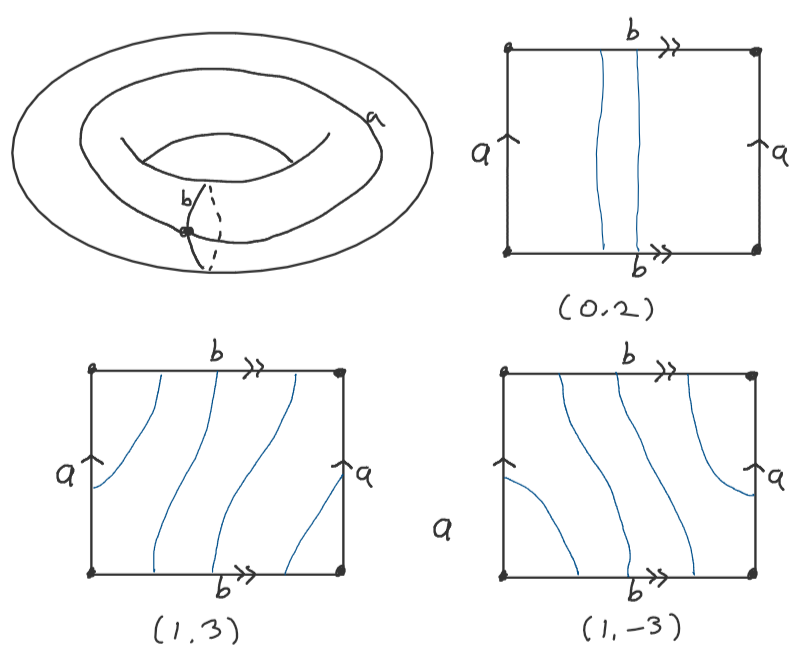}
    \caption{Some multicurves on $S_{1,1}$}
    \label{fig:ex3}
\end{figure}

From this we calculate $f_{a\cup b}$ as follows: $$f(x,y)=1+\frac{x}{1-x}+\frac{y}{1-y}+\frac{2xy}{(1-x)(1-y)}=\frac{1+xy}{(1-x)(1-y)}.$$Then $$f\left(x^{-1},y^{-1}\right)=\frac{1+x^{-1}y^{-1}}{\left(1-x^{-1}\right)\left(1-y^{-1}\right)}=\frac{1+xy}{(1-x)(1-y)}=f(x,y).$$
\end{example}
The theorem follows from the observation that multicurves can be viewed as integer points of a rational cone. In this case Stanley's reciprocity theorem on rational cones(\cite{MR411982}) is applicable. The following theorem, which can be found in \cite[p.~179, p.~188]{MR666158}, is one of its reformulations.

\begin{theorem}[Stanley]
\label{thm:Stanley}
For $r\leq m$, let $\Phi$ be an $r \times m$ integer matrix, $K$ be the set of its nonnegative integer solutions and $K^\circ$ its positive integer solutions. Then $k(z)=\sum\limits_{\alpha \in K}z^\alpha$, $k^\circ (z)=\sum\limits_{\alpha \in K^\circ}z^\alpha$ are rational functions. Factors of their denominators are of the form $\left(1-z^\beta\right)$ where $\beta$ is a nonnegative integer solution of $\Phi$. Furthermore, if $k^\circ \neq 0$, $k$ and $k^\circ$ satisfy  
$$k(z^{-1})=(-1)^{m-\mathrm{rank}(\Phi)}k^\circ(z).$$
\end{theorem}

\begin{remark}
Equivalently, we could define everything on $\Sigma_{g,n}$, the compact oriented surface of genus $g$ with $n$ boundaries, instead of $S_{g,n}$. In this case counting curve is defined to be a disjoint union of curves connecting between boundary curves just as in the case of ribbon graph. The number of curves counted will not be affected. In this point of view Theorem~\ref{thm:0} is a special case of Theorem~\ref{thm:main}, where the counting curves are defined to be those curves cutting across the handles. 
\end{remark}

\section{Proof}

\begin{lemma}
\label{lem1}
Suppose a counting curve $C$ is given. By adding additional curve components to $C$, we can obtain a new counting curve $C'$ such that $S_{g,n}$ is divided by open regions bounded by three (not necessarily distinct) curve components of $C'$ and those regions are homeomorphic to an open disk.
\end{lemma}
\textit{Proof.} Each open region bounded by $C$ is homeomorphic to an open disk and does not include any puncture. If not, it would allow non-contractible curves not intersecting $C$. Add diagonal curves to those polygon-shaped regions to divide them into smaller regions bounded by three curve components, and we have $C'.$ It is then obvious that $C'$ is a counting curve. \hspace*{\fill}$\square$

From now on, \emph{edge} refers to a component of $C'$ and \emph{triangle} refers to a region bounded by $C'$. It is possible that two of the three edges of a triangle coincide; see the triangle bounded by $c$ and $d$ in Figure~\ref{fig:ex2}. Since an edge bounds two triangles and a triangle is bounded by three edges (if counted properly), we have $2N$ triangles and $3N$ edges for some $N$. Then $2N$ triangles, $3N$ edges, $n$ points filling in the punctures consist a cell structure on a surface of genus $g$ and we have $2N-3N+n=2-2g$, hence $N\equiv n\ (mod\ 2)$.

Now observe how a multicurve appears in each triangles. When we consider a multicurve intersecting $C'$ at minimum in its isotopy class, their should be no bigons bounded by part of an edge and part of the multicurve. Thus such multicurve restricted to a triangle is a disjoint union of segments connecting between edges of the triangle. Its isotopy type is characterized by a triple of nonnegative integers representing the number of segments connecting each pair of edges. Since we have $2N$ triangles the configuration on the whole surface is represented by a $6N$-tuple of nonnegative integers. There must be same number of segments connected from each side of an edge so the tuple must satisfy linear relations corresponding to the edges. There are $3N$ edges and the tuple must be a nonnegative integer solution of some $3N \times 6N$ integer matrix $\Phi$. Conversely, if we have a nonnegative linear solution of such $\Phi$, we can glue those curves in triangles to obtain a multicurve. The argument above is suggests the following proposition.

\begin{proposition}
\label{prop1}
Given a counting curve $C'$ satisfying the condition of Lemma~\ref{lem1} with $2N$ triangles and $3N$ edges, there is a one-to-one correspondence between isotopy classes of multicurves and  nonnegative integer solutions of some $3N\times6N$ integer matrix $\Phi$. 
\end{proposition}

\begin{example}
In Figure~\ref{fig:ex2}, $C=a\cup b\cup c$ is a counting curve and can be refined to $C'=a\cup b\cup c\cup d\cup e\cup f$ which satisfies the condition of Lemma~\ref{lem1}. In the view of Proposition~\ref{prop1}, multicurves on $S_{1,2}$ corresponds(up to isotopy) to nonnegative integer solutions of the following $\Phi$. 
\begin{figure}[h]
    \centering
    \includegraphics[scale=0.75]{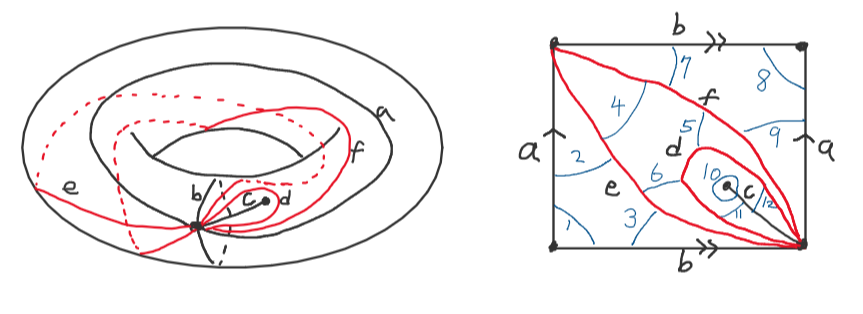}
    \caption{Edges and triangles on $S_{1,2}$}
    \label{fig:ex2}
\end{figure}
$$\Phi=\bordermatrix{%
&1&2&3&4&5&6&7&8&9&10&11&12\cr
a & 1 & 1 & 0 & 0 & 0 & 0 & 0 & -1 & -1 & 0 & 0 & 0\cr
b & 1 & 0 & 1 & 0 & 0 & 0 & -1 & -1 & 0 & 0 & 0 & 0\cr
c & 0 & 0 & 0 & 0 & 0 & 0 & 0 & 0 & 0 & 0 & 1 & -1\cr
d & 0 & 0 & 0 & 0 & 1 & 1 & 0 & 0 & 0 & 0 & -1 & -1\cr
e & 0 & 1 & 1 & -1 & 0 & -1 & 0 & 0 & 0 & 0 & 0 & 0\cr
f & 0 & 0 & 0 & 1 & 1 & 0 & -1 & 0 & -1 & 0 & 0 & 0\cr
}$$

\end{example}

\begin{figure}[h]
    \centering
    \includegraphics[scale=0.75]{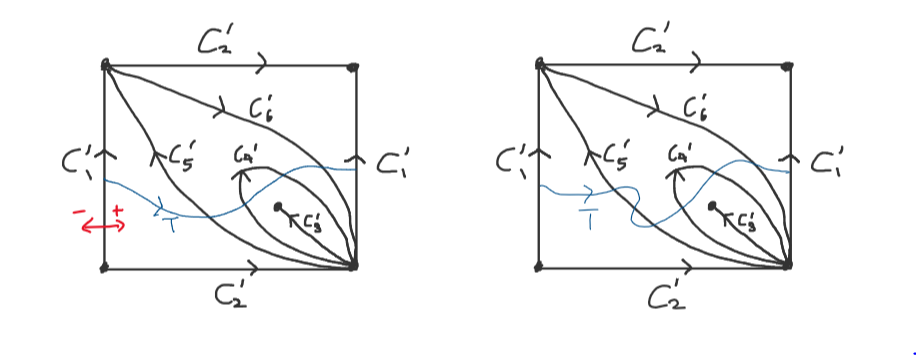}
    \caption{A curve configuration modified by an isotopy}
    \label{fig:ex4}
\end{figure}

\textit{Proof of Proposition~\ref{prop1}.} We only need to show that different solutions correspond to non-isotopic multicurves. Denote the edges by $C'_1, C'_2, \ldots C'_{3N}$. Fix their orientations. For an oriented connected curve component $T$ that meets $C'$ transversely finitely many times, we designate a cyclic word by moving along $T$ and concatenating free symbols $X_i$ or $X_i^{-1}$ (depending on the orientation) whenever it crosses $C_i$. Note that any isotopy of $T$ can be viewed as consecutively inserting or deleting $X_i X_i^{-1}$ or $X_i^{-1}X_i$'s(by making or removing bigons). For example, in Figure~\ref{fig:ex4} the cyclic word for the left figure is $\left[X_1 X_5 X_4 X_4^{-1} X_6\right]$ and after making a bigon(right figure) it becomes $\left[X_1 X_5 X_5^{-1} X_5 X_4 X_4^{-1} X_6\right]$. Hence the cyclic word, when viewed as a conjugacy class in the free group with generators $\left\{X_1, X_2, \ldots X_{3N}\right\}$, is well-defined up to isotopy. In the conjugacy class there is a cyclically reduced word which is unique up to cyclic shift(\cite[Theorem 1.3]{MR2109550}). This means $T$ is isotopic to a unique configuration without any bigon, hence no more than one solution of $\Phi$ can correspond to $T$. This applies to every multicurve by considering for each connected component.\hspace*{\fill}$\square$

From now on, fix $C$, $C'$, $\Phi$ as defined in $Proposition~\ref{prop1}$.

\begin{lemma}
\label{lem2}
$\mathrm{rank}(\Phi)=3N$.
\end{lemma}
\textit{Proof.} Let $v_1,v_2\cdots,v_{3N}$ be row vectors of $\Phi$. Suppose $\sum\limits_{i=1}^{3N} a_i v_i=0$. Choose a triangle. Pick columns corresponding to segments connecting between edges of the triangle. If the triangle has three distinct edges, say $v_1,v_2,v_3$. Then $\Phi$ looks as follows after reordering columns so that the picked columns come first, up to sign of rows.  
$$\left(\begin{array}{cccccc}
1&1&0&&\multirow{6}{*}{$\ast$}&\\
1&0&1&&&\\
0&1&1&&&\\
0&0&0&&&\\
&\vdots& &&&\\
0&0&0&&&
\end{array}\right)$$
From this, we have $a_1=a_2=a_3=0$.
Meanwhile, if two of edges coincide so that the chosen triangle is bounded by two distinct edges, say $v_1,v_2$, then $\Phi$ looks as follows up to reordering of rows and columns and sign of rows.

$$\left(\begin{array}{cccccc}
1&1&0&&\multirow{5}{*}{$\ast$}& \\
1&-1&0&& &\\
0&0&0&&&\\
&\vdots&&& &\\
0&0&0&& &
\end{array}\right)$$
Thus $a_1=a_2=0$. Since every edge bounds some triangle, we have $a_i=0$ for all $i$. We conclude that $v_1,v_2,\cdots,v_{3N}$ are independent and $\mathrm{rank}(\Phi)=3N$. \hspace*{\fill}$\square$

Next consider curves that are contractible to some puncture. On $S_{g,n}$ we have $n$ curves up to isotopy, say $d_1,d_2,\ldots, d_n$. Denote their corresponding $6N$-tuples by $\delta_1,\delta_2,\ldots,\delta_n$.

\begin{lemma}
\label{lem3}
$\sum\limits_{i=1}^n\delta_i=\left(1,1,\ldots,1\right)$.
\end{lemma}
\textit{Proof.} Consider a segment connecting between two (not necessarily distinct) edges of a triangle. It corresponds to an angle between those two edges. By collecting all segments we obtain a multicurve whose corresponding tuple is $\left(1,1,\ldots,1\right)$. This multicurve passes along each angle exactly once, so this is $\bigcup\limits_{i=1}^n d_i$ and the corresponding tuple is $\sum\limits_{i=1}^n\delta_i$. \hspace*{\fill}$\square$

\textit{Proof of Theorem~\ref{thm:main}.}

Note that $\left(1,1,\ldots,1\right)$ is the unique minimal positive integer solution of $\Phi$. By applying Theorem~\ref{thm:Stanley} we have
$$k\left(z^{-1}\right)=(-1)^{6N-\mathrm{rank}(\Phi)}k^\circ (z)=(-1)^{6N-\mathrm{rank}(\Phi)}z^{\left(1,1,\ldots,1\right)}k(z)$$
where $k(z)$ and $k^\circ(z)$ are defined as in Theorem~\ref{thm:Stanley}. Proposition~\ref{prop1} tells us that $g_{C'}(x)=k(z)$ under the substitution $z_{(i,j)}=x_i^{\frac{1}{2}} x_j^{\frac{1}{2}}$, where $x_i, x_j$ are variables corresponding to edges $C_i, C_j$ and $z_{(i,j)}$ is the variable corresponding to the segment connecting those edges. Define $\alpha_i$ by $x^{\alpha_i}=z^{\delta_i}$. Since every multicurve is a disjoint union of an essential multicurve and (possibly multiple copies of) $d_i$'s, we have $f_{C'}(x)\prod\limits_{i=1}^n\sum\limits_{k=0}^\infty x^{k\alpha_i}=g_{C'}(x)$, or equivalently $f_{C'}(x)=g_{C'}(x)\prod\limits_{i=1}^n\left(1-x^{\alpha_{i}}\right)$. With Lemma~\ref{lem2} and Lemma~\ref{lem3}, we have \begin{align*}
f_{C'}\left(x^{-1}\right)
&=g_{C'}\left(x^{-1}\right)\prod\limits_{i=1}^n\left(1-x^{-\alpha_{i}}\right)\\
&=(-1)^{6N-\mathrm{rank}(\Phi)}z^{\left(1,1,\ldots,1\right)}g_{C'}(x)\prod\limits_{i=1}^n\left(1-x^{-\alpha_i}\right)\\
&=(-1)^{N}\left(\prod\limits_{i=1}^nx^{\alpha_{i}}\right)g_{C'}(x)\prod\limits_{i=1}^n\left(1-x^{-\alpha_{i}}\right)\\
&=(-1)^n g_{C'}(x)\prod\limits_{i=1}^n\left(x^{\alpha_i}-1\right)\\
&=g_{C'}(x)\prod\limits_{i=1}^n\left(1-x^{\alpha_{i}}\right)\\&=f_{C'}(x).
\end{align*}
By substituting $1$ to variables corresponding to components of $C'$ that are not components of $C$, we obtain $g_C$, $f_C$ from $g_{C'}$, $f_{C'}$, respectively. This is justified by the following argument. The denominators of $g_{C'}$ are of the form $\left(1-x^\beta\right)$ by Theorem~\ref{thm:Stanley}. Since every nontrivial multicurve intersects $C$, at least one variable corresponding to a component of $C$ has a positive exponent in $x^\beta$. In other words, $\left(1-x^\beta\right)$ does not vanish by the substitution. Hence $g_C$ is well-defined. Since $f_{C'}(x)=g_{C'}(x)\prod\limits_{i=1}^n\left(1-x^{\alpha_{i}}\right)$, we can say the same thing about $f_{C}$. By the substitution on the equation $f_{C'}\left(x^{-1}\right)=f_{C'}(x)$, we conclude that $f_{C}\left(x^{-1}\right)=f_{C}(x)$. \hspace*{\fill}$\square$

\begin{acknowledgment}
The author would like to emphasize his sincerest gratitude towards his advisor Junho Peter Whang for introducing one of his valuable works \cite{MR4157427} and suggesting the interesting challenge of generalizing it. His numerous counsels during the writing of this paper must be also appreciated. 
\end{acknowledgment}

\renewcommand\refname{\vskip -1cm}
\section{References}
\bibliographystyle{unsrt}
\bibliography{cite}

\end{document}